\documentclass[11pt,a4paper]{article}

\usepackage{amsmath,amssymb}

\newcommand{\counte}{theorem}
\newtheorem{theorem} {\bf Theorem}[section]

\newtheorem{defn}[\counte]{\bf Definition}

\newtheorem{lemma}[\counte]{\bf Lemma}

\newtheorem{remark}[\counte]{\bf Remark}

\allowdisplaybreaks

\renewcommand{\thefootnote}{\fnsymbol{footnote}}

\begin{document}

\renewcommand{\thefootnote}{\arabic{footnote}}

\centerline{\Large\bf A Schur-Toponogov theorem in Riemannian geometry}
\centerline{\Large\bf  \& a new proof of Toponogov's theorem}
\centerline{\Large\bf in Alexandrov geometry}

\vskip5mm

\centerline{Yusheng Wang \footnote{Supported by NSFC 11471039.
\hfill{$\,$}}}

\vskip6mm

\noindent{\bf Abstract.} In the paper, we give a Schur-Toponogov theorem in Riemannian geometry,
which not only generalizes Schur's and Toponogov's theorem but also indicates their relation.
Inspired by its proof, we also supply a new proof of Toponogov's theorem  (in the large)
in Alexandrov geometry.

\vskip1mm


\vskip1mm


\vskip6mm

\setcounter{section}{-1}

\section {Introduction}

In Rimannian geometry, a kind of elegant geometry is on distance comparison by curvature;
for instance, Shcur's theorem ([Ho]). In the present paper, we denote by $[pq]$ a minimal
geodesic between $p$ and $q$ in a Riemannian manifold or  Alexandrov space.

\begin{theorem}[A. Shcur]\label{3.1}
Let $\gamma(s)|_{[0,L]}\subset \Bbb E^n$ ($n$-dimensional Euclidean space)
and $\tilde\gamma(s)|_{[0,L]}\subset \Bbb E^2$ be two arc-length parameterized $C^2$-curves,
and let $\kappa(s), \tilde\kappa(s)$ be their
absolute curvature respectively. If $\kappa(s)\leq\tilde\kappa(s)$ for all $s$,
and if $\tilde\gamma(s)|_{[0,L]}\cup[\tilde\gamma(0)\tilde\gamma(L)]$ is
a convex curve, then the distance $$|\gamma(0)\gamma(L)|\geq|\tilde\gamma(0)\tilde\gamma(L)|;$$
and if equality holds, then $\gamma(s)|_{[0,L]}$ is equal to $\tilde\gamma(s)|_{[0,L]}$\
up to an isometry of $\Bbb E^n$.
\end{theorem}

Recall that $\kappa(s)=|D_{\dot\gamma(s)}\dot\gamma(s)|,$ where $\dot\gamma(s)$ denotes the tangent vector of $\gamma(s)$ and
$D$ is the canonical connection of $\Bbb E^n$.  Theorem 0.1 also holds for piecewise $C^2$-curves ([Ho], cf. [Su]), i.e.,
$\gamma(s)$ and $\tilde\gamma(s)$ are not differential at only a finite number of points,
and at such a point `$\kappa(s)\leq\tilde\kappa(s)$' means that
 $|\dot{\gamma}_-(s)\dot{\gamma}_+(s)|\leq |\dot{\tilde\gamma}_-(s)\dot{\tilde\gamma}_+(s)|$ (where $\dot{\gamma}_\pm(s)$ denotes the right and left tangent vectors).

\vskip1mm

We now consider another significant result on distance comparison---Toponogov's
theorem ([BGP], [Pe], [GM]). In the paper, we always denote by $\Bbb S^n_k$ the complete and simply
connected $n$-dimensional Riemannian manifold with constant sectional curvature $k$.
For a given $[pq]$ in a Riemannian manifold (or Alexandrov space), we denote
by $\uparrow_p^q$ the direction from $p$ to $q$ of $[pq]$ (in a Riemannian manifold,
$\uparrow_p^q$ is just the unit tangent vector of $[pq]$ at $p$).

\begin{theorem}[Toponogov]\label{2.1}
Let $M^n$ be a complete Riemannian manifold with $\sec_M\geq k$. Let $p\in M$ and $[qr]\subset M$,
and let $\tilde p\in\Bbb S^2_k$ and $[\tilde q\tilde r]\subset \Bbb S^2_k$
such that $|\tilde p\tilde q|=|pq|$, $|\tilde q\tilde r|=|qr|$, and $|pq|+|pr|+|qr|<\frac{2\pi}{\sqrt k}$ for $k>0$.
Then the following holds:

\vskip1mm

\noindent{\rm (0.2.1)} If $|\tilde p\tilde
r|=|pr|$, then for $s\in[qr]$ and
$\tilde s\in[\tilde q\tilde r]$ with $|qs|=|\tilde q\tilde s|$,
$$|ps|\geq|\tilde p\tilde s|.\eqno{(0.1)}$$

\vskip1mm

\noindent{\rm (0.2.2)} If $|\uparrow_{q}^p\uparrow_{q}^r|=|\uparrow_{\tilde q}^{\tilde p}\uparrow_{\tilde q}^{\tilde r}|$ for some $[pq]$,
then for $s\in[qr]$ and
$\tilde s\in[\tilde q\tilde r]$ with $|qs|=|\tilde q\tilde s|$,
$$|ps|\leq|\tilde p\tilde s|.\eqno{(0.2)}$$

\vskip1mm

\noindent{\rm (0.2.3)}
If the equality in {\rm (0.1)} (resp. {\rm(0.2)}) holds for some $s\in [qr]^\circ$ \footnote{ $[qr]^\circ$ denotes the interior part of $[qr]$.}
(resp. $s\in [qr]\setminus\{q\}$), then there is $[pq]$ and $[pr]$ (resp. $[ps]$) which together with $[qr]$ (resp. $[pq]$ and $[qs]$) bounds a
convex surface which can be isometrically embedded into $\Bbb S^2_k$.
\end{theorem}

\begin{remark}\label{0.3} {\rm  (0.3.1) In Theorem 0.2, if $k>0$, then $|pq|+|pr|+|qr|\leq\frac{2\pi}{\sqrt k}$,
and equality holds if and only if $M^n$ is isometric to $\Bbb S^n_k$ with $p, q, r$ lying in a great circle (cf. [BGP]).

\vskip1mm

\noindent(0.3.2) It is not hard to see that (0.2.1) is equivalent to (0.2.2),
and (0.2.2) is equivalent to that, in (0.2.1),
$|\uparrow_{q}^p\uparrow_{q}^r|\geq|\uparrow_{\tilde q}^{\tilde p}\uparrow_{\tilde q}^{\tilde r}|$ for any $[pq]$.

\vskip1mm

\noindent(0.3.3) Theorem 0.2 has a corresponding version for $\sec_M\leq k$ with inverse
inequalities, but it requires that
$[qr]\subset M\setminus C(p)$, where $C(p)$ denotes the cut locus of $p$ (cf. [Me]).
If $[qr]\subset M\setminus C(p)$, the corresponding
`$|\uparrow_{q}^p\uparrow_{q}^r|\leq|\uparrow_{\tilde q}^{\tilde p}\uparrow_{\tilde q}^{\tilde r}|$' (see (0.3.2))
is an immediate corollary of Rauch's first comparison theorem (but this way does not work if $\sec_M\geq k$).}
\end{remark}

We now state our main result, a Schur-Toponogov theorem, which unifies and generalizes Theorem 0.1, 0.2 and (0.3.3).

\vskip2mm

\noindent {\bf Theorem A} {\it Let $M$ be a complete Riemannian manifold with $\sec_M\leq k$ (resp. $\sec_M\geq k$)
and let $\gamma(s)|_{[0,L]}\subset M$ and $\tilde \gamma(s)|_{[0,L]}\subset \Bbb S^2_k$
be two arc-length parameterized $C^2$-curves. Suppose that the absolute curvature of  $\gamma(s)$ and $\tilde\gamma(s)$,
$\kappa(s)$ and $\tilde \kappa(s)$, satisfy that  $\kappa(s)\leq\tilde \kappa(s)$
(resp. $\kappa(s)\geq\tilde \kappa(s)$) for all $s$. Then the following holds:

\vskip1mm

\noindent{\rm (A1)} If $\tilde\gamma(s)|_{[0,L]}\cup[\tilde\gamma(0)\tilde\gamma(L)]$
is a convex curve and $\gamma(s)|_{[0,L]}\subset M\setminus C(\gamma(0))$
(resp. if $\gamma(s)|_{[0,L]}$ is convex to $\gamma(0)$), then
$$|\gamma(0)\gamma(L)|\geq \text{(resp. $\leq$) } |\tilde \gamma(0)\tilde \gamma(L)|;\eqno{(0.3)}$$
and if equality holds and $|\gamma(0)\gamma(L)|+L<\frac{2\pi}{\sqrt k}$ for $k>0$, then
(resp. then there is $[\gamma(0)\gamma(s)]$ such that) $\bigcup_{s\in[0,L]}[\gamma(0)\gamma(s)]$
with induced metric is isometric to $\bigcup_{s\in[0,L]}[\tilde \gamma(0)\tilde \gamma(s)]$.

\vskip1mm

\noindent{\rm (A2)} Let $p\in M, \tilde p \in\Bbb S^2_k$ with $|p\gamma(0)|=|\tilde p\tilde\gamma(0)|$,
and let $|p\gamma(0)|+|p\gamma(L)|+L<\frac{2\pi}{\sqrt k}$ for $k>0$.
If $[\tilde p\tilde\gamma(0)]\cup\tilde\gamma(s)|_{[0,L]}\cup[\tilde p\tilde\gamma(L)]$
is a convex curve and $\gamma(s)|_{[0,L]}\subset M\setminus C(p)$ (resp. if $\gamma(s)|_{[0,L]}$ is convex to $p$),
then the following holds:

\vskip1mm

{\rm (A2.1)} If $|p\gamma(L)|=|\tilde p\tilde\gamma(L)|$, then for $s\in (0,L)$
$$|p\gamma(s)|\leq \text{(resp. $\geq$) } |\tilde p\tilde \gamma(s)|.\eqno{(0.4)}$$

{\rm (A2.2)} If (resp. If for some $[p\gamma(0)]$) $|\uparrow_{\gamma(0)}^p\dot\gamma(0)|=|\uparrow_{\tilde\gamma(0)}^{\tilde p}\dot{\tilde\gamma}(0)|$,
then for $s\in (0,L]$ $$|p\gamma(s)|\geq \text{(resp. $\leq$) }|\tilde p\tilde \gamma(s)|.\eqno{(0.5)}$$

{\rm (A2.3)} If the equality in {\rm (0.4)} or {\rm(0.5)} holds for some $s\in (0,L)$ or $s_0\in (0,L]$ respectively, then
(resp. then there is $[p\gamma(s)]$ such that) $\bigcup_{s\in[0,L]}[p\gamma(s)]$
or $\bigcup_{s\in[0,s_0]}[p\gamma(s)]$ with induced metric is isometric to $\bigcup_{s\in[0,L]}[\tilde p\tilde \gamma(s)]$ or
$\bigcup_{s\in[0,s_0]}[\tilde p\tilde \gamma(s)]$ respectively.}

\vskip3mm

In Theorem A, $\kappa(s)$ (the absolute curvature of  $\gamma(s)$) is still defined to be $|D_{\dot\gamma(s)}\dot\gamma(s)|$,
where $D$ is the Levi-Civita connection of $M$.
And the definition of a convex curve in $\Bbb S^2_k$
is the same as in $\Bbb E^2$, while the notion `$\gamma(s)|_{[0,L]}$ is convex to $p$ (or $\gamma(0)$)' will be
introduced in Section 2 (see Definition 2.1).

\begin{remark}\label{0.4} {\rm (0.4.1) From our proof, one can see that (A2.1) is equivalent to  (A2.2),
and (A2.2) is equivalent to that, in (A2.1),
$|\uparrow_{\gamma(0)}^p\dot\gamma(0)|\leq \text{ (resp. $\geq$) } |\uparrow_{\tilde\gamma(0)}^{\tilde p}\dot{\tilde\gamma}(0)|$
for any $[p\gamma(0)]$ (cf. (0.3.2)).

\vskip1mm

\noindent (0.4.2) From our proof, (A1) can be viewed as a corollary of (A2), but we formulate (A1) and (A2)
in a theorem because the proof of (A2.3) requires (0.3).}
\end{remark}

\begin{remark}\label{0.5} {\rm (0.5.1) It is clear that (A1) of Theorem A
includes Theorem 0.1. And (A2) for $\sec_M\geq k$ includes Theorem 0.2 (see Remark 2.4,
where a key point is that if $\gamma(s)|_{[0,L]}$ is a minimal geodesic
(so is $\tilde\gamma(s)|_{[0,L]}$), then it is convex to $p$ automatically).
We would like to point out that  Theorem
0.1 and 0.2 are not used in our proof; in other words, we supply a new proof for Theorem 0.1
and 0.2 in this paper.

\vskip1mm

\noindent (0.5.2) One can find counterexamples for Theorem A without the condition `$\gamma(s)|_{[0,L]}\subset M\setminus C(\gamma(0))$' or `$\gamma(s)|_{[0,L]}\subset M\setminus C(p)$' (resp. `$|\gamma(0)\gamma(L)|+L<\frac{2\pi}{\sqrt k}$'
or `$|p\gamma(0)|+|p\gamma(L)|+L<\frac{2\pi}{\sqrt k}$' for $k>0$) when $M$ is a cylinder
(resp. $\Bbb S^2_k$).

\vskip1mm

\noindent (0.5.3) Similar to Theorem 0.1, Theorem A also holds for piecewise $C^2$-curves
$\gamma(s)$ and $\tilde \gamma(s)$ (in fact, our proof works for such general case),
and at each non-differential point `$\kappa(s)\leq\ (\geq)\ \tilde\kappa(s)$' means
`$|\dot{\gamma}_-(s)\dot{\gamma}_+(s)|\leq\ (\geq)\  |\dot{\tilde\gamma}_-(s)\dot{\tilde\gamma}_+(s)|$'.}
\end{remark}

Toponogov's theorem  is not only a powerful tool in Riemannian geometry, but also
the base of Alexandrov geometry. An Alexandrov space $X$ with curvature $\geq k$
is roughly defined to be a locally compact length space
\footnote{If an Alexandrov space $X$ is locally compact,
then for any $p\in X$ there is a neighborhood $U$ of $p$ such that between any two points in $U$
there is a shortest path (geodesic). In [BGP],
the original definition of Alexandrov spaces is for locally complete spaces.} on which Theorem 0.2 without
the rigidity part (0.2.3) holds locally ([BGP]). In fact, if $X$ is complete, then Theorem 0.2
holds globally on $X$, for which without the rigidity part there are three proofs so far
([BGP], [Pl], [Sh]; refer to [GM] for a proof for the rigidity part). The proof in [BGP]
can be viewed as a version corresponding to the original definition of Alexandrov spaces with curvature $\geq k$
in [BGP] (which also adapts to locally complete spaces), while the proofs in [Pl] and [Sh] can be
viewed corresponding to (0.2.2) in Theorem 0.2.
Inspired by our proof of Theorem A (especially Lemma 1.2),
we supply a new proof which is a version corresponding to (0.2.1) in Theorem 0.2, i.e. we will prove:

\vskip2mm

\noindent {\bf Theorem B} {\it Let $X$ be a complete Alexandrov space with curvature $\geq k$.
Let $p\in X$ and $[qr]\subset X$,
and let $\tilde p\in\Bbb S^2_k$ and $[\tilde q\tilde r]\subset \Bbb S^2_k$
such that $|\tilde p\tilde q|=|pq|,|\tilde p\tilde
r|=|pr|$, $|\tilde q\tilde r|=|qr|$. Suppose that $|pq|+|pr|+|qr|<\frac{2\pi}{\sqrt k}$ for $k>0$.
Then for any $s\in[qr]$ and
$\tilde s\in[\tilde q\tilde r]$ with $|qs|=|\tilde q\tilde s|$, we have that $|ps|\geq|\tilde p\tilde s|.$}

\vskip2mm

Similar to (0.3.1), if $k>0$ in Theorem B, then one can show that $|pq|+|pr|+|qr|\leq\frac{2\pi}{\sqrt k}$
([BGP]).

\section{Proof of Theorem A for the case where $\sec_M\leq k$}

\subsection{Preparations}

\noindent{\bf a. Closed convex curves in $\Bbb S_k^2$}

Let $\tilde\gamma(s)|_{[0,L]}$ be an arc-length parameterized $C^2$-curve in $\Bbb S^2_k$, and $p\in \Bbb S^2_k$. If $[\tilde p\tilde\gamma(0)]\cup\tilde\gamma(s)|_{[0,L]}\cup[\tilde p\tilde\gamma(L)]$ is a convex curve, then for any $\tilde\gamma(s)$ and $\tilde q\in [\tilde p\tilde\gamma(0)]\cup\tilde\gamma(s)|_{[0,L]}\cup[\tilde p\tilde\gamma(L)]$,
$$\text{the angle between $D_{\dot{\tilde \gamma}(s)}\dot{\tilde \gamma}(s)$
and $[\tilde \gamma(s)\tilde q]$ at $\tilde \gamma(s)$ is not bigger than } \frac\pi2.\eqno{(1.1)}$$
Moreover,

\vskip2mm

\centerline{$|\tilde p\tilde\gamma(0)|+|\tilde p\tilde\gamma(L)|+L\leq \frac{2\pi}{\sqrt k}\text{ for }k>0$, and}
`$=$' holds iff $[\tilde p\tilde\gamma(0)]\cup\tilde\gamma(s)|_{[0,L]}\cup[\tilde p\tilde\gamma(L)]$ is the
union of two half great circles. \hfill(1.2)

\vskip2mm

\noindent In order to see (1.2), it suffices to show that a closed and convex piecewise minimal geodesic
$\Gamma\subset\Bbb S_k^2$ (i.e., a convex polygon)
is of length $\leq \frac{2\pi}{\sqrt k}$, and equality holds if and only if $\Gamma$ is the
union of two half great circles. Note that $\Gamma$ is the boundary of
the intersection of a finite number of half spheres.
Then by induction on the number of these half spheres, one can see the wanted property of $\Gamma$.

\vskip2mm

\noindent{\bf b. Index forms of Jacobi fields}

Let $M$ be a complete Riemannian manifold, and let $c(t)|_{[a,b]}\subset M$ be a normal geodesic.
Recall that for a normal \footnote{Here, `normal' means that $X(t)\perp \dot{c}(t)$} Jacobi field $X(t)$ along $c(t)$,
the index form
$$I_a^b(X,X)\triangleq\int_a^{b}\left(|\dot{X}(t)|^2-\langle R_{\dot c(t)X(t)}\dot c(t), X(t)\rangle\right)dt.$$

Let $p\in M$ and $\gamma(s)|_{(-\epsilon,\epsilon)}\subset M$ be
an arc-length parameterized $C^2$-curve. If $\gamma(0)\not\in C(p)$, then there is
a Jacobi field $U(t)$ along the normal geodesic $c(t)\triangleq[p\gamma(0)]$ with $U(0)=0$ and $U(|p\gamma(0)|)=\dot{\gamma}(0)$
such that
$$\frac{d|p\gamma(s)|}{ds}|_{s=0}=-\cos|\dot{\gamma}(0)\uparrow_{\gamma(0)}^p|$$
and
$$\frac{d^2|p\gamma(s)|}{d^2s}|_{s=0}=\left\langle D_{U(t)}U(t),\dot{c}(t)\right\rangle|_{0}^{|p\gamma(0)|}+
I_{0}^{|p\gamma(0)|}({U}^\perp,{U}^\perp),\eqno{(1.3)}$$
where $U^\perp(t)$ is the projection of $U(t)$ to the orthogonal space of $\dot{c}(t)$ in $T_{c(t)}M$.

\vskip2mm

\noindent{\bf c. Index form comparison}

On index forms, we have the following comparison result (cf. [Wu], [CE]), which plays an essential role in proving Rauch's first comparison theorem.

\begin{lemma}\label{2.1}
Let $M$ and $\tilde M$ be two complete $n$-dimensional Riemannian manifolds, and let $c(t)|_{[0,\ell]}\subset M$ and
$\tilde c(t)|_{[0,\ell]}\subset \tilde M$ be
two normal minimal geodesics. Let $K_{\max}(t)$ (resp. $\tilde K_{\min}(t)$) be the maximum (resp. the minimum) of sectional curvatures at $c(t)$ (resp. $\tilde c(t)$).  And let $J(t)$ and $\tilde J(t)$
be normal Jacobi fields along $c(t)$ and $\tilde c(t)$ respectively.
If $K_{\max}(t)\leq\tilde K_{\min}(t)$, $J(0)=\tilde J(0)=0$ and $|J(\ell)|=|\tilde J(\ell)|$,
then $$I_{0}^{\ell}(J,J)\geq I_{0}^{\ell}(\tilde J, \tilde J).$$
Moreover, if equality holds and if $\tilde M=\Bbb S^n_k$,
then $K(\dot c(t), J(t))\equiv k$ (where $K(\dot c(t), J(t))$ is the sectional curvature of the plane spanned by $\dot c(t)$ and $J(t)$).
\end{lemma}

\vskip2mm

\noindent{\bf d. Two easy facts on distance comparison}

\begin{lemma}\label{2.1}
Let $M$ and $\tilde M$ be two complete Riemannian manifolds, let $p\in M$ and $\tilde p\in \tilde M$, and let $c(t)|_{[0,\ell]}\subset M$ and $\tilde c(t)|_{[0,\ell]}\subset \tilde M$ be
two $c^1$-curves such that $|pc(0)|=|\tilde p\tilde c(0)|$ and $|pc(\ell)|=|\tilde p\tilde c(\ell)|$.
If $|pc(t)|-|\tilde p\tilde c(t)|$ attains a local minimum at $t_0\in (0,\ell)$
and there is a unique minimal geodesic between $\tilde p$ and $\tilde c(t_0)$, then
$$|\dot{c}(t_0)\uparrow_{c(t_0)}^{p}|=|\dot{\tilde c}(t_0)\uparrow_{\tilde c(t_0)}^{\tilde p}|.$$
\end{lemma}

\noindent{\it Proof}. This is an almost immediate corollary of the first variation formula.
\hfill $\Box$

\begin{lemma}\label{2.1}
Let $[p_1p_2]\subset\Bbb S^2_k$
with $|p_1p_2|<\frac{\pi}{\sqrt \kappa}$ for $\kappa>0$, and let $q\in [p_1p_2]^\circ$ and
$c(t)|_{[0,\ell]}\subset\Bbb S^2_k$ be a normal minimal geodesic with $c(0)=q$
and $\ell\leq\frac{\pi}{\sqrt \kappa}$ for $\kappa>0$. Then
$|p_1c(t)|+|p_2c(t)| \text{ is strictly increasing  with respect to } t$.
In particular, if $|qp_i|$ ($i=1$ or $2$) is sufficiently small compared to $|qc(t_0)|$,
then $|p_ic(t)|<|p_ic(t_0)|$ for all $t\in [0,t_0)$.
\end{lemma}

\noindent{\it Proof}. It is not hard to see that the lemma follows from the Law of Cosine.
\hfill $\Box$

\subsection {Proof of Theorem A for the case where $\sec_M\leq k$}

In this subsection, we always assume that $\sec_M\leq k$. And we remark that the following proof goes through for piecewise $C^2$-curves
$\gamma(s)$ and $\tilde \gamma(s)$ as mentioned in (0.5.3).

\vskip2mm

\noindent {\bf Step 1.} To prove (A2.1).

\vskip1mm

(This step is the essential part of the whole proof for Theorem A.)

\vskip1mm

Observe that `$|p\gamma(0)|+|p\gamma(L)|+L<\frac{2\pi}{\sqrt k}$' for $k>0$ implies that $|p\gamma (s)|<\frac{\pi}{\sqrt k}$ for all $s$.
It follows that for any $\bar s\in[0,L]$ there is a $C^2$-curve $\tilde \gamma_{\bar s}(s)|_{[0,L]}\subset\Bbb S_k^2$,
equal to $\tilde \gamma(s)|_{[0,L]}$ up to an isometry of $\Bbb S_k^2$, such that
$$|p\gamma(\bar s)|=|\tilde p\tilde\gamma_{\bar s}(\bar s)| \text{ and }
|\uparrow_{\gamma(\bar s)}^p\dot\gamma(\bar s)|=
|\uparrow_{\tilde\gamma_{\bar s}(\bar s)}^{\tilde p}\dot{\tilde\gamma}_{\bar s}(\bar s)|. \eqno{(1.4)}$$
Since $\gamma(s)|_{[0,L]}\subset M\setminus C(p)$,
$\gamma(s)|_{[0,L]}$ (and $\tilde \gamma_{\bar s}(s)|_{[0,L]}$ similarly) determines a Jacobi field $U_s(t)|_{[0,|p\gamma(s)|]}$ along
$\beta_s(t)|_{[0,|p\gamma(s)|]}\triangleq[p\gamma(s)]$ such that  $U_s(|p\gamma(s)|)=\dot\gamma(s)$.
Then putting (1.1), (1.4), Lemma 1.1 and (1.3) together, we can conclude that
$$\frac{d^2|p\gamma(s)|}{d^2s}\big|_{s=\bar s}\geq \frac{d^2|\tilde p\tilde \gamma_{\bar s}(s)|}{d^2s}\big|_{s=\bar s}\ ;\eqno{(1.5)}$$
and equality holds if and only if
$$\text{ either $|\uparrow_{\gamma(\bar s)}^p\dot\gamma(\bar s)|=0$ or $\pi$,}\eqno{(1.6)}$$
\centerline{or $\kappa(\bar s)=\tilde\kappa(\bar s)$ and the angle between $D_{\dot\gamma(\bar s)}\dot\gamma(\bar s)$
and $\uparrow_{\gamma(\bar s)}^p$ is equal to}

\hskip3cm that between $D_{\dot{\tilde\gamma}_{\bar s}(\bar s)}\dot{\tilde\gamma}_{\bar s}(\bar s)$
and $\uparrow_{\dot{\tilde\gamma}_{\bar s}(\bar s)}^{\tilde p}$.\hfill (1.7)

\vskip2mm

\noindent Moreover, when (1.7) occurs, $D_{\dot\gamma(\bar s)}\dot\gamma(\bar s)$ lies in the plane
spanned by $\dot\gamma(\bar s)$ and $\uparrow_{\gamma(\bar s)}^p$ in $T_{\gamma(\bar s)}M$,
and the sectional curvature of the plane spanned by
$\dot\beta_{\bar s}(t)$ and $U_{\bar s}(t)$ in $T_{\beta_{\bar s}(t)}M$
$$K(\bar s,t)=\kappa \text{ for all } t\in[0,|p\gamma(\bar s)|]\eqno{(1.8)}.$$

Then we will finish Step 1 through the following three substeps.

\vskip2mm

\noindent{\bf Substep 1.} To give a proof for the case where $[0,L]$ is one kind of ``best''
intervals, i.e.: the `$\geq$' in (1.5) is `$>$' for all $\bar s$, and
for any $s_1\neq s_2\in [0,L]$ there is
a $C^2$-curve $\tilde \gamma_{s_1,s_2}(s)|_{[0,L]}\subset\Bbb S_k^2$,
which is equal to $\tilde \gamma(s)|_{[0,L]}$ up to an isometry of $\Bbb S_k^2$, such that
$$|p\gamma(s_i)|=|\tilde p\tilde\gamma_{s_1,s_2}(s_i)|\eqno{(1.9)}$$
and
$$[\tilde p\tilde\gamma_{s_1,s_2}(s_1)]\cup\tilde\gamma_{s_1,s_2}(s)|_{[s_1,s_2]}\cup[\tilde p\tilde\gamma_{s_1,s_2}(s_2)]
 \text{ is convex}.\eqno{(1.10)}$$
(For example, $[0,L]$ will be such a ``best'' interval if $\tilde \gamma(s)|_{[0,L]}$
(and so $\gamma(s)|_{[0,L]}$) is a minimal geodesic.)

We now assume that (A2.1) is not true for such a ``best'' case.
Our strategy is to find some $[a,b]\subset[0,L]$ with $b-a$ sufficiently small such that
$$ |\uparrow_{\gamma(a)}^p\dot\gamma(a)|>
|\uparrow_{\tilde\gamma_{a,b}(a)}^{\tilde p}\dot{\tilde\gamma}_{a,b}(a)|
\text{ and }
\frac{d^2|p\gamma(s)|}{d^2s}\big|_{[a,b]}>\frac{d^2|\tilde p\tilde \gamma_{a,b}(s)|}{d^2s}\big|_{[a,b]},\eqno{(1.11)}$$
which contradicts to `$|p\gamma(b)|=|\tilde p\tilde\gamma_{a,b}(b)|$' (see (1.9)).

By the assumption right above, we can let $s_0\in (0,L)$ such that
$$|p\gamma(s_0)|-|\tilde p\tilde \gamma(s_0)|=
\max_{s\in(0,L)}\{|p\gamma(s)|-|\tilde p\tilde \gamma(s)|\}>0.$$
By Lemma 1.2,
$$|\uparrow_{\gamma(s_0)}^p\dot\gamma(s_0)|=
|\uparrow_{\tilde\gamma(s_0)}^{\tilde p}\dot{\tilde\gamma}(s_0)|.\eqno{(1.12)}$$
Note that we can prolong $[\tilde\gamma(s_0)\tilde p]$ to
$[\tilde\gamma(s_0)\tilde p']$ such that
$|\tilde\gamma(s_0)\tilde p'|=|\gamma(s_0)p|$.
Note that $[\tilde\gamma(0)\tilde\gamma(L)]\cap[\tilde\gamma(s_0)\tilde p]\neq\emptyset$ due to the convexity
of $[\tilde p\tilde\gamma(0)]\cup\tilde\gamma(s)|_{[0,L]}\cup[\tilde p\tilde\gamma(L)]$.
Then by Lemma 1.3, at least one of `$|\tilde p'\tilde\gamma(0)|>|p\gamma(0)|$' and
`$|\tilde p'\tilde\gamma(L)|>|p\gamma(L)|$' holds;
moreover,

\vskip1mm

\centerline {$|\tilde p'\tilde\gamma(0)|>|p\gamma(0)|$ (resp. $|\tilde p'\tilde\gamma(L)|>|p\gamma(L)|$) if $s_0$ (resp. $L-s_0$) }
\hskip 5cm is sufficiently small.\hfill(1.13)

\vskip1mm

\noindent We now assume that $|\tilde p'\tilde\gamma(L)|>|p\gamma(L)|$, and consider the curve
$\tilde\gamma_{s_0,L}(s)|_{[0,L]}$ (due to that $[0,L]$ is ``best''). `$|\tilde p'\tilde\gamma(L)|>|p\gamma(L)|$' implies that
$|\uparrow_{\tilde\gamma(s_0)}^{\tilde p}\dot{\tilde\gamma}(s_0)|>
|\uparrow_{\tilde\gamma_{s_0,L}(s_0)}^{\tilde p}\dot{\tilde\gamma}_{s_0,L}(s_0)|$, so by (1.12)
$$|\uparrow_{\gamma(s_0)}^p\dot\gamma(s_0)|>
|\uparrow_{\tilde\gamma_{s_0,L}(s_0)}^{\tilde p}\dot{\tilde\gamma}_{s_0,L}(s_0)|.\eqno{(1.14)}$$

By (1.14) and the first variation formula, it is easy to see that there is $s_0'\in (s_0,L)$ such that
$|p\gamma(s_0')|-|\tilde p\tilde\gamma_{s_0,L}(s_0')|=
\max_{s\in(s_0,L)}\{|p\gamma(s)|-|\tilde p\tilde \gamma_{s_0,L}(s)|\}>0.$
And consequently, we can set $[s_1,s_1']\triangleq [s_0,s_0']$ or $[s_0',L]$ such that
$$ |\uparrow_{\gamma(s_1)}^p\dot\gamma(s_1)|>
|\uparrow_{\tilde\gamma_{s_1,s_1'}(s_1)}^{\tilde p}\dot{\tilde\gamma}_{s_1,s_1'}(s_1)| \text{ or }
|\uparrow_{\gamma(s_1')}^p\dot\gamma(s_1')|<
|\uparrow_{\tilde\gamma_{s_1,s_1'}(s_1')}^{\tilde p}\dot{\tilde\gamma}_{s_1,s_1'}(s_1')| .
$$
By repeating the above process, we can get $[s_1, s_1']\supsetneq\cdots\supsetneq[s_i,s_i']\supsetneq\cdots$
such that
$$ |\uparrow_{\gamma(s_i)}^p\dot\gamma(s_i)|>
|\uparrow_{\tilde\gamma_{s_i,s_i'}(s_i)}^{\tilde p}\dot{\tilde\gamma}_{s_i,s_i'}(s_i)| \text{ or }
|\uparrow_{\gamma(s_i')}^p\dot\gamma(s_i')|<
|\uparrow_{\tilde\gamma_{s_i,s_i'}(s_i')}^{\tilde p}\dot{\tilde\gamma}_{s_i,s_i'}(s_i')| .
\eqno{(1.15)}$$
Moreover, by (1.13)
we can select $s_i$ and $s_i'$ such that they converge to some $\bar s$ as $i\to\infty$, and then
we put $[a,b]\triangleq[s_i,s_i']$ for a sufficiently large $i$.

We now need to check that $[a,b]$ satisfies (1.11). Note that by (1.15)
we can have the first inequality of (1.11). As for the second one,
since we have assumed that the `$\geq$' in (1.5) is `$>$', it suffices to show that
$$\tilde \gamma_{s_i,s_i'}(s)|_{[0,L]} \text{ converges to } \tilde \gamma_{\bar s}(s)|_{[0,L]}
\text{ as $i\to\infty$} \eqno{(1.16)}$$
(note that $\tilde \gamma_{s_i,s_i'}(s)|_{[0,L]}$ and $\tilde \gamma_{\bar s}(s)|_{[0,L]}$
are all equal to $\tilde \gamma(s)|_{[0,L]}$ up to isometries of $\Bbb S^2_k$).
Observe that the smaller $|s_i-s_i'|$ is the smaller $\left||\uparrow_{\gamma(s_i)}^p\dot\gamma(s_i)|-
|\uparrow_{\tilde\gamma_{s_i,s_i'}(s_i)}^{\tilde p}\dot{\tilde\gamma}_{s_i,s_i'}(s_i)|\right|$
is because both $\gamma(s)$ and $\tilde \gamma(s)$ are $C^2$-curves, and thus (1.16) follows.

\vskip2mm

\noindent{\bf Substep 2.} To give a proof for the case where $[0,L]$ is another ``best''
interval, namely, the `$\geq$' in (1.5) is `$=$' for all $\bar s$ and (1.6) does not occur at any $\bar s$.

Since $\gamma(s)|_{[0,L]}\subset M\setminus C(p)$ and $|p\gamma (s)||_{[0,L]}<\frac{\pi}{\sqrt k}$,
we can draw a $C^2$-curve $\bar{\tilde \gamma}(s)|_{[0,L]}$ in $\Bbb S_k^2$ with
$s$ being arc-length parameter such that
$|p\gamma(s)|=|\tilde p\bar{\tilde \gamma}(s)|$ for all $s\in[0,L]$.
Consider the corresponding (1.3) to $|\tilde p\bar{\tilde \gamma}(s)|$ together with (1.7) and (1.8),
and notice that (1.6) does not occur at any $\bar s$ in this substep.
It has to hold that $\bar{\tilde \gamma}(s)|_{[0,L]}$ is equal to $\tilde \gamma(s)|_{[0,L]}$
up to an isometry of $\Bbb S_k^2$, and thus
$|p\gamma(s)|=|\tilde p\tilde \gamma(s)|$ for all $s\in[0,L]$.

\vskip2mm

\noindent{\bf Substep 3.} To give a proof for general cases.

We first observe that (A2.1) holds if (1.6) occurs at all $\bar s\in [0,L]$ (which implies that (1.5) is an equality).
In fact, in such a situation, $\gamma(s)|_{[0,L]}$ has to be a minimal geodesic with $[p\gamma(1)]=[p\gamma(0)]\cup\gamma(s)|_{[0,L]}$ or $[p\gamma(0)]=[p\gamma(1)]\cup\gamma(s)|_{[0,L]}$, and $\tilde \gamma(s)|_{[0,L]}$ has the same phenomenon.

We now assume that $W\triangleq\{\bar s\in [0,L]|\ \text{(1.6) occurs at } \bar s\}$
is not equal to $[0,L]$.  Note that if (A2.1) is not true,
in order to get a contradiction we can assume that
$$|p\gamma(s)|>|\tilde p\tilde\gamma(s)| \text{ for all } s\in (0,L),\eqno{(1.17)}$$ which implies that
$$|\uparrow_{\gamma(0)}^p\dot\gamma(0)|\geq |\uparrow_{\tilde\gamma(0)}^{\tilde p} \dot{\tilde\gamma}(0)|\text{ (and } |\uparrow_{\gamma(L)}^p\dot\gamma(L)|\leq |\uparrow_{\tilde\gamma(L)}^{\tilde p} \dot{\tilde\gamma}(L)|).\eqno{(1.18)}$$

In the following, we will derive a contradiction under (1.18) (with (1.17)).

Note that
$W$ is a closed
subset of $[0,L]$, so $W^c\triangleq[0,L]\setminus W$ is a union of some intervals.
Due to the $C^2$-property of $\gamma(s)|_{[0,L]}$ and $\tilde\gamma(s)|_{[0,L]}$,
for any $\bar s\in W^c$ (i.e., $|\uparrow_{\gamma(\bar s)}^p\dot\gamma(\bar s)|\neq0$ or $\pi$),
there is an interval $I\subset W^c$ containing $\bar s$ such that for any  $s_{1},s_{2}\in I$ there exists $\tilde \gamma_{s_1,s_2}(s)|_{[0,L]}$
(equal to $\tilde \gamma(s)|_{[0,L]}$ up to an isometry of $\Bbb S_k^2$) satisfying (1.9) and (1.10).
Then by Substep 1 and 2, $W^c$ must contain a {\it nice} interval
$(s_1,s_2)$. Here, we call $(s_1,s_2)$ a nice interval if $\tilde \gamma_{s_1,s_2}(s)|_{[0,L]}$
satisfying (1.9) and (1.10) exists and $|p\gamma(s)|\leq |\tilde p\tilde \gamma_{s_1,s_2}(s)|$ for any $s\in (s_1,s_2)$.

We now can take a maximal nice interval $(s_{11},s_{21})\subset (0,L)$ (i.e. any $(c,d)\supsetneq (s_{11},s_{21})$ is not nice), which may contain the point in $W$.
It is clear that
$$|\uparrow_{\gamma(s_{11})}^p\dot\gamma(s_{11})|\leq |\uparrow_{\tilde\gamma_{s_{11},s_{21}}(s_{11})}^{\tilde p} \dot{\tilde\gamma}_{s_{11},s_{21}}(s_{11})|,\ |\uparrow_{\gamma(s_{21})}^p\dot\gamma(s_{21})|\geq |\uparrow_{\tilde\gamma_{s_{11},s_{21}}(s_{21})}^{\tilde p} \dot{\tilde\gamma}_{s_{11},s_{21}}(s_{21})|.\eqno{(1.19)}$$
By (1.17), $[s_{11}, s_{21}]\neq [0, L]$, so outside of $[s_{11}, s_{21}]$ it is possible that
there exists a maximal nice interval. Note that there are at most $\{(s_{1j},s_{2j})\}_{j=1}^\infty$
such that $(s_{1(j+1)},s_{2(j+1)})$ is a maximal nice interval
outside of $\bigcup_{l=1}^j[s_{1j}, s_{2j}]$ and
$$W\cup\bigcup_{j}[s_{1j},s_{2j}]=[0, L], \text{ and any $s\not\in W$ lies in some $[s_{1j}, s_{2j})$ and $(s_{1j'}, s_{2j'}]$}.\eqno{(1.20)}$$
In order to complete Substep 3, we introduce the following function
$$d(s)\triangleq\begin{cases} |\tilde p\tilde\gamma_{s_{1j},s_{2j}}(L)|, & s\in (s_{1j},s_{2j})\\
|\tilde p\tilde\gamma_{s}(L)|, & s\in [0,L]\setminus\bigcup_{j}(s_{1j},s_{2j}) \end{cases}$$
(for $\tilde\gamma_s$ refer to the beginning of Step 1).
Note that $d(L)=|\tilde p\tilde\gamma(L)|$, and (1.18) implies that
$d(0)\geq|\tilde p\tilde\gamma(L)|\ (=d(L))$.
However, we claim that
$$\text{ $d(s)\geq d(0)$ for all $s$ with $d(L)>d(0)$}, \eqno{(1.21)}$$
a contradiction (and thus Substep 3 is finished).

We now need only to verify (1.21)
(in fact, the proof implies that $d(s)|_{[0,L]}$ is non-decreasing).
We first show that $d(s)\geq d(0)$ for all $s$.
Let $A$ be the collection of such $\bar s\in [0,L]$
that  $d(s)\geq d(0)$ on $[0,\bar s]$, and that for $s'\in(s_{1j},s_{2j})\cap [0,\bar s]$
(resp. $s'\in[0,\bar s]\setminus\bigcup_{j}(s_{1j},s_{2j})$) $[\tilde p\tilde\gamma_{s_{1j},s_{2j}}(s')]\cup \tilde\gamma_{s_{1j},s_{2j}}(s)|_{[s',L]}\cup[\tilde p\tilde\gamma_{s_{1j},s_{2j}}(L)]$
(resp. $[\tilde p\tilde\gamma_{s'}(s')]\cup \tilde\gamma_{s'}(s)|_{[s',L]}\cup[\tilde p\tilde\gamma_{s'}(L)]$)  is convex.
In fact, we can prove that $A$ is a non-empty, open and closed subset of $[0,L]$,
i.e. $A=[0,L]$.
Note that $A\neq\emptyset$ because $0\in A$ by (1.18). As for `open and closed',
it suffices to show that
$$\text{ if $\bar s\in A$ and $\bar s\neq L$, then there is a $\delta>0$ such that $[\bar s, \bar s+\delta)\subset A$} \eqno{(1.22)}$$
and
$$\text{ if $(\bar s-\delta', \bar s)\subset A$ for some $\delta'>0$, then $\bar s\in A$}. \eqno{(1.23)}$$

Case 1: If $\bar s$ lies in some $(s_{1j},s_{2j})$, then (1.22) and (1.23) hold automatically.

Case 2: If $\bar s$ is equal to some $s_{1j}$ (resp. $s_{2j}$), then the first  (resp. the second)
inequality of corresponding (1.19) (together with the convexity of $\tilde\gamma(s)|_{[\bar s,L]}$) implies (1.22) (resp. (1.23)).
Moreover, if $\bar s=s_{1j}=s_{2j'}$, then the maximality of each $(s_{1j},s_{2j})$ implies that at least one
inequality of corresponding (1.19) is strict, and thus
$$d(s)\geq d(\bar s)>d(s') \text{ or } d(s)>d(\bar s)\geq d(s'),\ \forall\ s\in(s_{1j},s_{2j}), s'\in(s_{1j'},s_{2j'}). \eqno{(1.24)}$$
(Note that if $[0,L]=\bigcup_{j}[s_{1j},s_{2j}]$, then  Case 1 and 2 right above implies (1.21).)

\vskip1mm

By (1.20), the residual case is:

Case 3: For (1.22) (resp. (1.23)), $\bar s\in W\setminus\bigcup_{j}(s_{1j},s_{2j})$ and $\bar s\neq s_{1j}$
(resp. $\bar s\neq s_{2j}$) for any $j$. In this case, due to the similarity,
we just give a proof for (1.22).
We first observe that $|\uparrow_{\gamma(\bar s)}^p\dot\gamma(\bar s)|=\pi$.
Otherwise, $|\uparrow_{\gamma(\bar s)}^p\dot\gamma(\bar s)|=0$ (see (1.6)), so the convexity of
$[\tilde p\tilde\gamma_{\bar s}(\bar s)]\cup \tilde\gamma_{\bar s}(s)|_{[\bar s,L]}\cup[\tilde p\tilde\gamma_{\bar s}(L)]$
(note that $\bar s\in A\setminus \bigcup_{j}(s_{1j},s_{2j})$) implies that $\tilde\gamma_{\bar s}(s)|_{[\bar s,L]}$ is a minimal geodesic
($\subseteq [\tilde\gamma_{\bar s}(\bar s)\tilde p]$), so is $\gamma(s)|_{[\bar s,L]}$ ($\subseteq [\gamma(\bar s)p]$);
which implies that $\bar s$ has to be $L$ (otherwise, $(\bar s, L)$ belongs to some $(s_{1j},s_{2j})$).
Due to $|\uparrow_{\gamma(\bar s)}^p\dot\gamma(\bar s)|=\pi$, it is not hard to see that
there is a small $\delta>0$ such that one of the following holds:

\noindent Subcase 1: $[\bar s, \bar s+\delta)\subset W\setminus\bigcup_{j}(s_{1j},s_{2j})$;

\noindent Subcase 2: For any $\delta'\leq\delta$, $(\bar s, \bar s+\delta')$ contains an infinite number of $(s_{1j},s_{2j})$.

In Subcase 1, $\gamma(s)|_{[\bar s,\bar s+\delta)}$
is a minimal geodesic (belonging to $[p\gamma(\bar s+\delta)]$), but $\tilde\gamma_{\bar s}(s)|_{[\bar s,\bar s+\delta')}$ with $\delta'\leq\delta$
is not (otherwise, $(\bar s, \bar s+\delta')$ belongs to some $(s_{1j},s_{2j})$). Then  for any $\hat s\in (\bar s, \bar s+\delta)$, we can rotate
$\tilde\gamma_{\hat s}(s)|_{[0,L]}$ around $\tilde\gamma_{\hat s}(\hat s)$ to $\tilde\gamma^*_{\hat s}(s)|_{[0,L]}$ such that
$\tilde\gamma^*_{\hat s}(\bar s)\in [\tilde p\tilde\gamma_{\hat s}(\hat s)]$.
Note that $|\tilde p\tilde\gamma^*_{\hat s}(\bar s)|>|\tilde p\tilde\gamma_{\bar s}(\bar s)|$.
By comparing $\tilde\gamma_{\hat s}(s)|_{[0,L]}$, $\tilde\gamma^*_{\hat s}(s)|_{[0,L]}$ and $\tilde\gamma_{\bar s}(s)|_{[0,L]}$
(which are all isometric to the convex $\tilde\gamma(s)|_{[0,L]}$), it is not hard to see that
$$d(\hat s)=|\tilde p\tilde\gamma_{\hat s}(L)|>|\tilde p\tilde\gamma^*_{\hat s}(L)|>|\tilde p\tilde\gamma_{\bar s}(L)|=d(\bar s)\eqno{(1.25)}$$
and $[\tilde p\tilde\gamma_{\hat s}(\hat s)]\cup \tilde\gamma_{\hat s}(s)|_{[\hat s,L]}\cup[\tilde p\tilde\gamma_{\hat s}(L)]$ is convex.
It follows that $[\bar s, \bar s+\delta) \subset A$.

In Subcase 2, an ideal model is that $(\bar s, \bar s+\delta)$
belongs to a union of many $[s_{1j},s_{2j}]$. In this case, by a similar reason for Case 1 and 2,
we can conclude that $[\bar s, \bar s+\delta) \subset A$ and (1.24) implies that
$$d(\bar s+\delta)>d(\bar s).\eqno{(1.26)}$$
In other cases, similarly, we only need to show that  $\hat s\in A$,
where $\hat s\in (\bar s, \bar s+\delta)\cap W\setminus\bigcup_{j}(s_{1j},s_{2j})$ and $\hat s\neq s_{2j}$ for any $j$.
Moreover, we can assume that
$$\text{for any $\hat s'\in [\bar s,\hat s)$, $(\hat s',\hat s)$ is not a nice interval}; \eqno{(1.27)}$$
otherwise, since $\bar s,\hat s\not\in\bigcup_{j}(s_{1j},s_{2j}) $, we can reselect the maximal nice intervals in $(\bar s,\hat s)$ such that $\hat s$
is just some $s_{2j}$.
Note that when $\delta$ is sufficiently small,  $|\uparrow_{\gamma(\hat s)}^p\dot\gamma(\hat s)|=\pi$
because $|\uparrow_{\gamma(\bar s)}^p\dot\gamma(\bar s)|=\pi$ and $\hat s\in W$.
Similarly, we can rotate
$\tilde\gamma_{\hat s}(s)|_{[0,L]}$ around $\tilde\gamma_{\hat s}(\hat s)$ to $\tilde\gamma^*_{\hat s}(s)|_{[0,L]}$ such that
$\tilde\gamma^*_{\hat s}(\bar s)\in [\tilde p\tilde\gamma_{\hat s}(\hat s)]$.
Then similar to Subcase 1, we can conclude that $[\bar s, \bar s+\delta) \subset A$
and the corresponding (1.25) holds if we can show that $|\tilde p\tilde\gamma^*_{\hat s}(\bar s)|>|\tilde p\tilde\gamma_{\bar s}(\bar s)|$.
In fact, if $|\tilde p\tilde\gamma^*_{\hat s}(\bar s)|\leq |\tilde p\tilde\gamma_{\bar s}(\bar s)|$,
then $\tilde\gamma^*_{\hat s}(s)|_{[0,L]}$ has to be equal to some $\tilde\gamma_{\hat s',\hat s}(s)|_{[0,L]}$
with $\hat s'\in [\bar s,\hat s)$ and $(\hat s',\hat s)$ is a nice interval,
which contradicts (1.27).

Note that we have proven  $d(s)|_{[0,L]}\geq d(0)$, and the proof (especially (1.24-26)) implies  $d(L)>d(0)$,
i.e. (1.21) holds (so Substep 3 (and thus Step 1) is completed).

\vskip2mm

\noindent {\bf Step 2.} To prove (A2.2), i.e. $|p\gamma(s)|\geq|\tilde p\tilde \gamma(s)|\ \forall\ s\in (0,L]$;
and if equality holds for some $s_0\in (0,L]$, then $|p\gamma(s)|=|\tilde p\tilde \gamma(s)|$ for all $s\in [0,s_0]$.

We first give an observation: there is a $\delta>0$ such that $|p\gamma(s)|\geq|\tilde p\tilde \gamma(s)|$
for all $s\in (0,\delta)$.
In fact, due to $|\uparrow_{\gamma(0)}^p\dot\gamma(0)|=|\uparrow_{\tilde\gamma(0)}^{\tilde p}\dot{\tilde\gamma}(0)|$
(the condition in (A2.2)), for sufficiently small $\delta>0$ if $|p\gamma(\delta)|<|\tilde p\tilde \gamma(\delta)|$, then $\tilde \gamma_{0,\delta}(s)|_{[0,L]}$
satisfying (1.9) and (1.10) exists and
$|\uparrow_{\gamma(0)}^p\dot\gamma(0)|>|\uparrow_{\tilde\gamma_{0,\delta}(0)}^{\tilde p}\dot{\tilde\gamma}_{0,\delta}(0)|$.
However, by applying (A2.1) on $\gamma(s)|_{[0,\delta]}$ and $\tilde \gamma_{0,\delta}(s)|_{[0,\delta]}$, it has to hold that
$|\uparrow_{\gamma(0)}^p\dot\gamma(0)|\leq|\uparrow_{\tilde\gamma_{0,\delta}(0)}^{\tilde p}\dot{\tilde\gamma}_{0,\delta}(0)|$,
a contradiction.

Due to the observation right above, we can let $s_0\neq 0$ be the maximal $s$ such that  $|p\gamma(s)|\geq|\tilde p\tilde \gamma(s)|$
for all $s\in [0,s_0]$. In order to prove (A2.2), we need to show that $s_0=L$. If it is not true,
then $|p\gamma(s_0)|=|\tilde p\tilde \gamma(s_0)|$.
Then by applying (A2.1) on $\gamma(s)|_{[0,s_0]}$ and $\tilde \gamma(s)|_{[0,s_0]}$,
we have that $|p\gamma(s)|\leq|\tilde p\tilde \gamma(s)|$ for all $s\in [0,s_0]$.
It follows that $|p\gamma(s)|=|\tilde p\tilde \gamma(s)|$ for all $s\in [0,s_0]$,
which implies that  $|\uparrow_{\gamma(s_0)}^p\dot\gamma(s_0)|=|\uparrow_{\tilde\gamma(s_0)}^{\tilde p}\dot{\tilde\gamma}(s_0)|$.
Thereby, similar to the observation right above, there is a $\delta>0$ such that $|p\gamma(s)|\geq|\tilde p\tilde \gamma(s)|$
for all $s\in (s_0,s_0+\delta)$, which contradicts the maximality of $s_0$ (i.e., $s_0=L$).

Note that the proof for `$s_0=L$' implies that if $|p\gamma(s_0)|=|\tilde p\tilde \gamma(s_0)|$
with $s_0\in (0,L]$, then $|p\gamma(s)|=|\tilde p\tilde \gamma(s)|$ for all $s\in [0,s_0]$.

\vskip2mm

\noindent {\bf Step 3.} To prove (0.3) in (A1), i.e. $|\gamma(0)\gamma(L)|\geq|\tilde \gamma(0)\tilde \gamma(L)|$.

We will derive a contradiction by assuming $|\gamma(0)\gamma(L)|<|\tilde\gamma(0) \tilde\gamma(L)|$.

We first claim that $|\gamma(0)\gamma(s)|<|\tilde\gamma(0) \tilde\gamma(s)|$ for all $s\in (0,L)$.
If it is not true, then there is $s_0\in (0,L)$ such that $|\gamma(0)\gamma(s_0)|=|\tilde\gamma(0) \tilde\gamma(s_0)|$.
Note that if $k>0$, by (1.2), `$|\gamma(0)\gamma(L)|<|\tilde\gamma(0) \tilde\gamma(L)|$' implies that
$|\gamma(0)\gamma(L)|+L<\frac{2\pi}{\sqrt k}$, and thus
$|\gamma(0)\gamma(s_1)|+|\gamma(0)\gamma(s_1)|+|s_2-s_1|<\frac{2\pi}{\sqrt k}$ for all  $s_1, s_2\in [0,L]$.
Hence, we can apply (A2.1) on $\gamma(s)|_{[0,s_0]}$ and $\tilde\gamma(s)|_{[0,s_0]}$
(with $p=\gamma(0)$ and $\tilde p=\tilde \gamma(0)$) to
conclude that $|\gamma(0)\gamma(s)|\leq|\tilde\gamma(0) \tilde\gamma(s)|$ for all $s\in (0,s_0)$.
In fact, this implies that $|\gamma(0)\gamma(s)|\leq|\tilde\gamma(0) \tilde\gamma(s)|$ for all $s\in (0,L)$.
Hence, by Lemma 1.2 it holds that
$|\uparrow_{\gamma(s_0)}^{\gamma(0)}\dot\gamma(s_0)|=|\uparrow_{\tilde\gamma(s_0)}^{\tilde \gamma(0)}\dot{\tilde\gamma}(s_0)|$ .
Then we can apply (A2.2) on $\gamma(s)|_{[s_0,L]}$ and $\tilde\gamma(s)|_{[s_0,L]}$
(with $p=\gamma(0)$ and $\tilde p=\tilde \gamma(0)$) to conclude that $|\gamma(0)\gamma(L)|\geq|\tilde\gamma(0)\tilde \gamma(L)|$,
which contradicts `$|\gamma(0)\gamma(L)|<|\tilde\gamma(0) \tilde\gamma(L)|$'.

By the claim right above, we can select a sufficiently small $s_1\in(0,L)$ such that
$|\gamma(0)\gamma(s_1)|-|\tilde\gamma(0) \tilde\gamma(s_1)|=\min_{s\in [0,s_1]}\{
|\gamma(0)\gamma(s)|-|\tilde\gamma(0) \tilde\gamma(s)|\}$, which implies that
$$\left|\uparrow_{\gamma(s_1)}^{\gamma(0)}(-\dot\gamma(s_1))\right|\geq \left|\uparrow_{\tilde\gamma(s_1)}^{\tilde\gamma(0)}(-\dot{\tilde \gamma}(s_1))\right|.\eqno{(1.28)}$$
On the other hand, since $s_1$ is sufficiently small, there is $q\in M$ such that
$$|q\gamma(s_1)|=|\gamma(0)\gamma(s_1)| \text{ and } |q\gamma(0)|+|q\gamma(s_1)|=|\tilde\gamma(0)\tilde\gamma(s_1)|,\eqno{(1.29)}$$
and that $|\uparrow_{\gamma(s_1)}^q(-\dot\gamma(s_1))|=|\uparrow_{\gamma(s_1)}^q\uparrow_{\gamma(s_1)}^{\gamma(0)}|+|\uparrow_{\gamma(s_1)}^{\gamma(0)}(-\dot\gamma(s_1))|$,
which together with (1.28) implies that
$$\left|\uparrow_{\gamma(s_1)}^q(-\dot\gamma(s_1))\right|> \left|\uparrow_{\tilde\gamma(s_1)}^{\tilde\gamma(0)}(-\dot{\tilde \gamma}(s_1))\right|.\eqno{(1.30)}$$
Let $\tilde q\in[\tilde\gamma(0)\tilde\gamma(s_1)]$ such that $|\tilde q\tilde \gamma(s_1)|=|q\gamma(s_1)|\ (=|\gamma(0)\gamma(s_1)|)$,
which enables us to apply (A2.1) here to conclude that $|q\gamma(s)|\leq|\tilde q\tilde\gamma(s)|$ for all  $s\in [0,s_1]$.
However, by (1.30) it is clear that $|q\gamma(s)|>|\tilde q\tilde\gamma(s)|$ for $s$ less than and close to $s_1$, a contradiction.

\vskip2mm

\noindent {\bf Step 4.} To prove (A2.3) and the rigidity part of (A1).

We first prove the rigidity part for (A2.1), i.e. if $|p\gamma(s_0)|=|\tilde p\tilde\gamma(s_0)|$ for some $s_0\in(0,L)$ in (A2.1), then $\bigcup_{s\in[0,L]}[p\gamma(s)]$ with induced metric is isometric to $\bigcup_{s\in[0,L]}[\tilde p\tilde \gamma(s)]$.
Note that the proof will be done once one show that for any sufficiently small interval $[s_1,s_2]\subset[0,L]$,
$\bigcup_{s\in[s_1,s_2]}[p\gamma(s)]$ is convex and is isometric to $\bigcup_{s\in[s_1,s_2]}[\tilde p\tilde \gamma(s)]$.

Note that if $[s_1,s_2]$ is sufficiently small, we can assume that $\bigcup_{s\in[s_1,s_2]}[p\gamma(s)]$
lies in $M\setminus C(q)$ for any $q\in \bigcup_{s\in[s_1,s_2]}[p\gamma(s)]$. Then it suffices to verify that
$$|a_1a_2|=|\tilde a_1\tilde a_2| \text{ for any } a_i\in [p\gamma(s_i)] \text{ and }
\tilde a_i\in [\tilde p\tilde\gamma(s_i)] \text{ with } |pa_i|=|\tilde p\tilde a_i|, \eqno{(1.31)}$$
and
$$[a_1a_2] \subset \bigcup_{s\in[s_1,s_2]}[p\gamma(s)]. \eqno{(1.32)}$$

We first show that (1.31) is true if $a_2=\gamma(s_2)$, i.e.
$$|a_1\gamma(s_2)|=|\tilde a_1\tilde \gamma(s_2)|.\eqno{(1.33)}$$
From Step 1, we know that $|p\gamma(s)|\leq|\tilde p\tilde\gamma(s)|$ for all $s\in(0,L)$, so by Lemma 1.2
$$\left|\uparrow_{\gamma(s_0)}^p\dot\gamma(s_0)\right|=\left|
\uparrow_{\tilde\gamma(s_0)}^{\tilde p}\dot{\tilde\gamma}(s_0)\right|.$$
Then by the rigidity part of Step 2, we have that
$$|p\gamma(s)|=|\tilde p\tilde\gamma(s)| \text{ for all } s\in[0,L],\eqno{(1.34)}$$
which in turn implies that
$$\left|\uparrow_{\gamma(s)}^p\dot\gamma(s)\right|=\left|
\uparrow_{\tilde\gamma(s)}^{\tilde p}\dot{\tilde\gamma}(s)\right| \text{ for all } s\in [0,L]. \eqno{(1.35)}$$
By (1.34) and (1.35), we can apply (A2.1) on $\{\text{the curve } \gamma(s)|_{[s_1,s_2]}\cup[p\gamma(s_1)], \gamma(s_2)\}$
and $\{\text{the curve } \tilde\gamma(s)|_{[s_1,s_2]}\cup[\tilde p\tilde \gamma(s_1)],\tilde\gamma(s_2)\}$
to conclude that
$|\gamma(s_1)\gamma(s_2)|\leq |\tilde\gamma(s_1)\tilde \gamma(s_2)|$.
On the other hand, by Step 3 we know that
$|\gamma(s_1)\gamma(s_2)|\geq |\tilde\gamma(s_1)\tilde\gamma(s_2)|$, and thus $|\gamma(s_1)\gamma(s_2)|=|\tilde\gamma(s_1)\tilde\gamma(s_2)|$.
Then (1.33) follows (just as (1.34)), by which we can similarly apply (A2.1) on $\{\text{the curve } \gamma(s)|_{[s_1,s_2]}\cup[p\gamma(s_2)], a_1\}$
and $\{\text{the curve } \tilde\gamma(s)|_{[s_1,s_2]}\cup[\tilde p\tilde\gamma(s_2)], \tilde a_1\}$ to conclude (1.31).

For (1.32), we let $q\in[a_1a_2]$ and $\tilde q\in [\tilde a_1\tilde a_2]$ with $|qa_1|=|\tilde q\tilde a_1|$,
and assume that $\tilde q\in[\tilde p\tilde\gamma(\bar s)]$ with $\bar s\in[s_1,s_2]$.
Note that (1.31) implies that $|\uparrow_{a_1}^{\gamma(s_1)}\uparrow_{a_1}^{a_2}|=|\uparrow_{\tilde a_1}^{\tilde\gamma(s_1)}\uparrow_{\tilde a_1}^{\tilde a_2}|$.
This together with (1.35) enables us to apply the rigidity part of Step 2 on
$\{\text{the curve } [a_1\gamma(s_1)]\cup\gamma(s)|_{[s_1,s_2]}\cup[a_2\gamma(s_2)], q\}$ and
$\{\text{the curve } [\tilde a_1\tilde \gamma(s_1)]\cup\tilde\gamma(s)|_{[s_1,s_2]}\cup[\tilde a_2\tilde\gamma(s_2)], \tilde q\}$,
and $\{\text{the curve } [a_1p]\cup[a_2p], q\}$ and $\{\text{the curve } [\tilde a_1\tilde p]\cup[\tilde a_2\tilde p], \tilde q\}$
respectively to conclude that $q$ lies in $[p\gamma(\bar s)]$. That is, (1.32) follows.

\vskip1mm

The rigidity part for (A2.2) can be seen similarly because the corresponding (1.34) and (1.35) to $[0,s_0]$ hold
(note that by the rigidity part of Step 2, if the equality in (0.5) holds for some $s_0\in (0,L]$,
then $|p\gamma(s)|=|\tilde p\tilde \gamma(s)|$ for all $s\in [0,s_0]$).

\vskip1mm

Similarly, the rigidity part of (A1) can be verified once we show that $|\gamma(0)\gamma(s)|=|\tilde \gamma(0)\tilde\gamma(s)|$ for all $s\in [0,L]$.
In fact, due to `$|\gamma(0)\gamma(L)|=|\tilde \gamma(0)\tilde\gamma(L)|$', we can apply (A2.1)
on $\{\text{the curve }\gamma(s)|_{[0,L]}, \gamma(0)\}$ and $\{\text{the curve } \tilde\gamma(s)|_{[0,L]},\tilde\gamma(0)\}$
to conclude that $|\gamma(0)\gamma(s)|\leq|\tilde \gamma(0)\tilde\gamma(s)|$ for all $s\in [0,L]$.
On the other hand, by Step 3 we can see that $|\gamma(0)\gamma(s)|\geq|\tilde \gamma(0)\tilde\gamma(s)|$,
so it follows that $|\gamma(0)\gamma(s)|=|\tilde \gamma(0)\tilde\gamma(s)|$.
\hfill$\Box$

\section{Proof of Theorem A for the case where $\sec_M\geq k$}

In this section, we always assume that $\sec_M\geq k$.

We first introduce the notion `$\gamma(s)|_{[0,L]}$ is convex to $p$' in Theorem A.
Let $\gamma(s)|_{[0,L]}\subset M$ be an arc-length parameterized $C^2$-curve. Note that given $p\in M$, $|p\gamma(s)||_{[0,L]}$
is a 1-Lipschitz function (maybe not a $C^1$-function). Then for $k\leq 0$ and if
$|p\gamma(s)||_{[0,L]}<\frac{\pi}{\sqrt k}$ for $k>0$, given $\bar p\in \Bbb S_k^2$, there is a unique (up to a rotation)
arc-length parameterized Lipschitz-curve $\bar{\gamma}(s)|_{[0,L]}$ in $\Bbb S_k^2$
such that
$$\text{$|p\gamma(s)|=|\bar p\bar{\gamma}(s)|$ for all $s\in[0,L]$},\eqno {(2.1)}$$
and $[\bar p\bar\gamma(s)]$ turns clockwise as $s$ increases ([PP], in which
$\bar{\gamma}(s)|_{[0,L]}$ is called a development of $\gamma(s)|_{[0,L]}$ with respect to $p$).
If $k>0$ and $|p\gamma(s)|=\frac{2\pi}{\sqrt k}$ for some $s\in [0,L]$, by the Maximum Diameter Theorem ([CE]),
we know that $M$ is isometric to $\Bbb S_k^n$;
and thus there is an arc-length parameterized $C^2$-curve $\bar{\gamma}(s)|_{[0,L]}$ such that (2.1) holds.

\begin{defn}\label{2.1} {\rm The $\gamma(s)|_{[0,L]}$ right above is said to be
{\it convex to $p$} if

\vskip1mm

\noindent (2.1.1) for all $s\in [0,L]$, there is
$[p\gamma(s)]$ such that the angle between $D_{\dot\gamma(s)}\dot\gamma(s)$ and $\uparrow_{\gamma(s)}^p$
is less than or equal to $\frac\pi2$, and $\{\dot\gamma(s), D_{\dot\gamma(s)}\dot\gamma(s),\uparrow_{\gamma(s)}^p\}$
lies in a plane of $T_{\gamma(s)}M$ \footnote{If $M=\Bbb S^n_k$, (2.1.1) implies that
$p$ and $\gamma(s)|_{[0,L]}$ fall in an $\Bbb S^2_k\subset\Bbb S^n_k$.}.

\vskip1mm

\noindent(2.1.2) $[\bar p\bar\gamma(0)]\cup\bar\gamma(s)|_{[0,L]}\cup[\bar p\bar\gamma(L)]$ is a convex curve in $\Bbb S_k^2$.}
\end{defn}

\begin{remark}\label{2.2} {\rm  (2.2.1) In Definition 2.1, if $\gamma(s)|_{[0,L]}$ is a piecewise $C^2$-curve, then
at non-differential points, $\dot\gamma(s)$ shall be replaced by $\dot\gamma_\pm(s)$ and it shall be added that
$$|(-\dot\gamma_-(s))\uparrow_{\gamma(s)}^p|+|\dot\gamma_+(s)\uparrow_{\gamma(s)}^p|=|(-\dot\gamma_-(s))\dot\gamma_+(s)|.$$
(2.2.2) If $k>0$ and $|p\gamma(s)|=\frac{2\pi}{\sqrt k}$ for some $s\in [0,L]$ (so $M$ is isometric to $\Bbb S_k^n$),
then by (1.2), (2.1.2) implies that each of $[p\gamma(0)]\cup\gamma(s)|_{[0,L]}\cup[p\gamma(L)]$ and
$[\bar p\bar\gamma(0)]\cup\bar\gamma(s)|_{[0,L]}\cup[\bar p\bar\gamma(L)]$
is a union of two half great circles. }
\end{remark}

Since Theorem A for $\sec_M\geq k$ is almost an immediate corollary of the proof in Section 1,
we only supply a rough proof for it.

\vskip2mm

\noindent{\bf Proof of Theorem A for $\sec_M\geq k$}.

Since $\gamma(s)|_{[0,L]}$ is convex to $\gamma(0)$ and $p$ in (A1) and (A2),
for the $\bar\gamma(s)|_{[0,L]}$ and $\bar p$ corresponding to $\gamma(s)|_{[0,L]}$ and $p$,
$\bar\gamma(s)|_{[0,L]}\cup[\bar \gamma(0)\bar\gamma(L)]$ and
$[\bar p\bar\gamma(0)]\cup\bar\gamma(s)|_{[0,L]}\cup[\bar p\bar\gamma(L)]$ are convex curves in $\Bbb S_k^2$ (by (2.1.2)).
We observe that the proof of Theorem A for $\sec_M\leq k$ in Section 1  applies to $\tilde \gamma(s)|_{[0,L]}$ and
$\bar \gamma(s)|_{[0,L]}$ (which both lie in $\Bbb S_k^2$)
by replacing $\gamma(s)|_{[0,L]}$ and $\tilde \gamma(s)|_{[0,L]}$ (in Section 1) by
$\tilde \gamma(s)|_{[0,L]}$ and $\bar \gamma(s)|_{[0,L]}$ respectively.
A problem here is that $\bar\gamma(s)|_{[0,L]}$ is only a 1-Lipschitz curve. However,
via (1.3), (2.1.1) implies that the corresponding (1.5) still holds (in the support sense),
and then one can check that the whole proof in Section 1 still works here.
This together with (2.1) implies that (0.3-5) follows. Moreover, if equalities hold in (0.3-5), then
$\tilde \gamma(s)|_{[0,L]}$ has to be equal to $\bar \gamma(s)|_{[0,L]}$
up to an isometry of $\Bbb S_k^2$, which implies $\bar \gamma(s)|_{[0,L]}$ is also
a $C^2$-curve. It then is not hard to see the rigidity part of (A1) and (A2.3).
\hfill$\Box$

\vskip2mm

We will end this section by showing that Definition 2.1 is not so artificial
through the following two remarks (especially Remark 2.4) on it.

\begin{remark}\label{2.3} {\rm  (2.1.1) implies that $\bar\gamma(s)\subset\Bbb S_k^2$ satisfies that
$$ \text{$\forall\ s\in [0,L]$, $\exists\ \delta>0$ s.t.
$[\bar p\bar\gamma(s-\delta)]\cup\bar\gamma(s)|_{[s-\delta,s+\delta]}\cup[\bar p\bar\gamma(s+\delta)]$ is convex
}\eqno{(2.2)}$$
(here if $s=0$ (resp. $L$), then $s-\delta$ (resp. $s+\delta$) should be $0$  (resp. $L$)).
Recall that for $p\in M$ ($\sec_M\geq k$) and an arc-length parameterized minimal geodesic $\sigma(s)|_{[a,b]}\subset M$,
if we let $f(s)=\rho_k(|p\sigma(s)|)$ with $\rho_k(x)=\begin{cases} \frac{1}{k}(1-\cos(\sqrt kx)), & k>0\\
\frac{x^2}{2}, & k=0\\
\frac{1}{k}(1-\cosh(\sqrt{-k}x)), & k<0
\end{cases}$, then
$$f''(s)\leq 1-kf(s) \text{ for all } s\in[a,b]\ (\text{in the support sense, [Pe]}).\eqno{(2.3)}$$
Thereby, for each $s_0\in[0,L]$, by comparing $\gamma(s)|_{[0,L]}$ and the minimal geodesic tangent to $\gamma(s)|_{[0,L]}$
at $\gamma(s_0)$,  via (1.3) one can see that (2.1.1) implies that (2.3) holds for $\bar p$ and $\bar\gamma(s)|_{[0,L]}$.
On the other hand, according to [PP], (2.2) is equivalent to (2.3) for $\bar p$ and $\bar\gamma(s)|_{[0,L]}$ ($\subset\Bbb S_k^2$).
It then follows that (2.1.1) implies (2.2).}
\end{remark}

\begin{remark}\label{2.4} {\rm  (2.4.1) If $\gamma(s)|_{[0,L]}$ is a minimal geodesic in Definition 2.1,
then $\gamma(s)|_{[0,L]}$ is convex to $p$ automatically. It is clear that $\gamma(s)|_{[0,L]}$ satisfies (2.1.1), so
it suffices to verify (2.1.2). We only need to consider the case where
$\gamma(0)\not\in[p\gamma(1)]$ and $\gamma(1)\not\in [p\gamma(0)]$.
In such a situation, for the $\bar p$ and $\bar\gamma(s)|_{[0,L]}$
corresponding to $p$ and $\gamma(s)|_{[0,L]}$, $s<|\bar p\bar\gamma(0)|+|\bar p\bar\gamma(s)|$ and thus
$|\uparrow_{\bar p}^{\bar\gamma(0)}\uparrow_{\bar p}^{\bar\gamma(s)}|<\pi$ for all $s\in [0,L]$.
This together with (2.2) (see Remark 2.3) implies (2.1.2).

\vskip1mm

\noindent(2.4.2) In (A2) of Theorem A for the case where $\sec_M\geq k$, if $\gamma(s)|_{[0,L]}$ is a minimal geodesic,
then $\tilde \gamma(s)|_{[0,L]}$ has to be a minimal geodesic too by `$\kappa(s)\geq\tilde\kappa(s)$'.
Then due to (2.4.1), (A2.1) and (A2.2) imply (0.2.1) and (0.2.2) respectively. Moreover,
note that $[p\gamma(0)]$ and $[p\gamma(L)]$ are also convex to $\gamma(L)$
and $\gamma(0)$ respectively. Hence, if the equality of (0.4) holds for some $s\in(0,L)$,
then by (A2.3) and by applying (A2) of Theorem A on $\{[p\gamma(0)], \gamma(L)\}$ and $\{[\tilde p\tilde \gamma(0)], \tilde \gamma(L)\}$
(and  $\{[p\gamma(L)], \gamma(0)\}$ and $\{[\tilde p\tilde \gamma(L)], \tilde \gamma(0)\}$),
it is not hard to conclude (0.2.3). That is, (A2) of Theorem A includes Theorem 0.2.}
\end{remark}

\section{Proof of Theorem B}

\subsection{Preparations}

Note that a proof of Theorem B has to depend only on its local version, i.e. the definition of
Alexandrov spaces with lower curvature bound, which enables us to define the angle between two minimal geodesics
starting from a common point ([BGP]).

In this section, $X$ always denotes a complete Alexandrov space with curvature $\geq k$.
Given $[pq],[pr]\subset X$, we denote by $\angle qpr$ the angle
between them (i.e. $|\uparrow_p^q\uparrow_p^r|$). And we
denote by $\triangle qpr$, a triangle, the union of three minimal geodesics
$[pq], [pr], [qr]\subset X$.

From the definition of angles, one can see the following easy properties ([BGP]).

\vskip2mm

\noindent {\bf Lemma 3.1} {\it Let $[pq],[rr']\subset X$  with $r\in[pq]^\circ$. Then
$\angle prr'+\angle qrr'=\pi$.}

\vskip2mm

\noindent {\bf Lemma 3.2} {\it Let  $\{[p_iq_i]\}_{i=1}^\infty$ and $\{[p_ir_i]\}_{i=1}^\infty$
be two sequences of minimal geodesics in $X$.
If $[p_iq_i]\to [pq]$ and $[p_ir_i]\to [pr]$ as $i\to\infty$, then $\angle qpr\leq \liminf\limits_{i\to\infty}\angle q_ip_ir_i$.}

\vskip2mm

Furthermore, we have the following easy observation.

\vskip2mm

\noindent {\bf Lemma 3.3} {\it Let $\{[p_iq_i]\}_{i=1}^\infty,\{[p_ir_i]\}_{i=1}^\infty\subset X$,
and let $p\in X$. If $|p_ip|\to 0$, $|p_ir_i|\to 0$ and $|p_iq_i|\to a>0$ as $i\to\infty$,
then we have that
$$|q_ir_i|\leq |q_ip_i|-|p_ir_i|\cos\angle q_ip_ir_i+o(|p_ir_i|).\eqno{(3.1)}$$}
\noindent{\it Proof}. By Toponogov's theorem around $p$,
there is $\bar q_i\in [p_iq_i]$ near $p$ such that
$$|\bar q_ir_i|\leq|\bar q_ip_i|-|p_ir_i|\cos\angle q_ip_ir_i+o(|p_ir_i|).$$
Then (3.1) follows from that $|q_ir_i|\leq |q_i\bar q_i|+|\bar q_ir_i|$.
\hfill $\Box$

\vskip2mm

In Lemma 3.3, if $p_i=p$, $q_i=q$, and $[p_ir_i]$ lie in a minimal geodesic,
then we can select $[pq]$ such that (3.1) is an equality (the first variation formula, [BGP]).
We will end this subsection by Alxandrov's lemma ([BGP]), a basic tool in Alexandrov geometry.

\vskip2mm

\noindent {\bf Lemma 3.4} {\it Let $\triangle pqr$, $\triangle pqs, \triangle abc \subset\Bbb S^2_k$,
where $\triangle pqr$ and $\triangle pqs$ are joined to each other in an exterior way along $[pq]$, such that
$|ab|=|pr|$, $|ac|=|ps|$, $|bc|=|qr|+|qs|$, and $|ab|+|ac|+|bc|<\frac{2\pi}{\sqrt k}$ if $k>0$. Then
$\angle pqr+\angle pqs\leq \pi$ (resp. $\geq\pi$) if and only if  $\angle prq\geq\angle abc$ and  $\angle psq\geq\angle acb$
(resp. $\angle prq\leq\angle abc$ and  $\angle psq\leq\angle acb$).}

\subsection{Proof of Theorem B}

Due to the similarity, we only give a proof for the case where $k=0$.

Assume that Theorem B is not true, i.e. there is $s\in [qr]$
and $\tilde s\in[\tilde q\tilde r]$ with $|qs|=|\tilde q\tilde s|$ such that
$|ps|<|\tilde p\tilde s|.$ To get a contradiction, roughly, we will construct a sufficiently small triangle
which does not satisfy Toponogov's theorem. The process
is completed through the following three steps, where for $\triangle abc\subset X$ we denote by $\tilde\angle abc$
the corresponding angle of its comparison triangle
\footnote{We usually call $\triangle \tilde a\tilde b\tilde c\subset\Bbb S^2_k$
the comparison triangle of $\triangle abc$ if $|\tilde a\tilde b|=|ab|, |\tilde a\tilde c|=|ac|$ and $|\tilde b\tilde c|=|bc|$.}  in $\Bbb S^2_k$.

\vskip1mm

\noindent{\bf Step 1}. To find $r'\in[qr]^\circ$ such that,
for any $[pr']$, $\angle pr'q<\tilde \angle pr'q$ or $\angle pr'r<\tilde \angle pr'r$.

\vskip1mm

By Step 1, there is a $\triangle p_1q_1r_1$ (with $p_1, q_1, r_1\in \triangle pqr$) such that
$$\angle p_1q_1r_1<\tilde\angle p_1q_1r_1.\eqno{(3.2)}$$
Let
$\text{peri}(p_1q_1r_1)$ and $(p_1q_1r_1)_{\min}$ denote the perimeter (i.e. $|p_1q_1|+|p_1r_1|+|q_1r_1|$) and
the length of the minimal side of $\triangle p_1q_1r_1$ respectively.

\vskip1mm

\noindent{\bf Step 2}. To find a $\triangle p_2q_2r_2$ with $p_2, q_2, r_2\in \triangle p_1q_1r_1$ such that
$\text{peri}(p_2q_2r_2)<\text{peri}(p_1q_1r_1)$,
$\angle p_2q_2r_2<\tilde \angle p_2q_2r_2$, and either

\noindent (3.3) $\cos\angle p_2q_2r_2-\cos\tilde \angle p_2q_2r_2\geq\cos\angle p_1q_1r_1-\cos\tilde \angle p_1q_1r_1$, or

\noindent (3.4) the following holds:

(3.4.1) $\text{peri}(p_2q_2r_2)<\text{peri}(p_1q_1r_1)-\frac12(p_1q_1r_1)_{\min}$,

(3.4.2) if $(p_2q_2r_2)_{\min}<(p_1q_1r_1)_{\min}$, then $\text{peri}(p_2q_2r_2)<4.5(p_1q_1r_1)_{\min}$.

\vskip1mm

\noindent{\bf Step 3}. To repeat Step 2 to get $\{\triangle p_iq_ir_i\}_{i=1}^\infty$
each of which satisfies $\text{peri}(p_{i}q_{i}r_{i})<\text{peri}(p_{i-1}q_{i-1}r_{i-1})$,
$\angle p_{i}q_{i}r_{i}<\tilde \angle p_{i}q_{i}r_{i}$, and the corresponding (3.3) or (3.4).

\vskip1mm

Note that any subsequence of $\{\triangle p_iq_ir_i\}_{i=1}^\infty$ contains a converging subsequence if $\{\triangle p_iq_ir_i\}_{i=1}^\infty$ lie in a compact domain of $X$.
It is obviously true if $X$ is compact. Anyway, due to our method of locating $\triangle p_iq_ir_i$,
it turns out to be true even when $X$ is not compact (see Remark 3.5 below).

Thereby, if $\text{peri}(p_iq_ir_i)\to 0$ as $i\to\infty$, i.e. $\triangle p_iq_ir_i$ converge to a point
passing to a subsequence,
then `$\angle p_iq_ir_i<\tilde \angle p_iq_ir_i$' contradicts the local Toponogov's theorem (cf. (0.3.2)).
In fact, either all $(p_iq_ir_i)_{\min}$ have a positive lower bound,
or $(p_iq_ir_i)_{\min}\to 0$ passing to a subsequence;
and thus if there is an infinite number of $\triangle p_{i}q_{i}r_{i}$ satisfying (3.4),
it is easy to see that $\text{peri}(p_iq_ir_i)\to 0$.

Then we can assume that $\text{peri}(p_iq_ir_i)$ converges to a positive number as $i\to\infty$,
and meanwhile the corresponding (3.3) occurs for all $i\geq 2$.
In this case, passing to a subsequence, $\triangle p_iq_ir_i$ converge to a minimal geodesic
or a triangle as $i\to\infty$. For the former case,
it is not hard to see that (3.3) for each $i$ together with Lemma 3.3 implies a contradiction.
For the latter case,  it is clear that (3.3) for each $i$ and Lemma 3.2 imply that
$\triangle p_iq_ir_i$ converge to a triangle $\triangle\bar p\bar q\bar r$ with
$$\cos\angle \bar p\bar q\bar r-\cos\tilde \angle \bar p\bar q\bar r\geq\cos\angle p_1q_1r_1-\cos\tilde \angle p_1q_1r_1;$$
and thus we can repeat such a process on $\triangle\bar p\bar q\bar r$ until a contradiction is gotten.

\vskip1mm

In the rest of the proof, we need only to show how we accomplish Step 1 and 2.

\vskip1mm

\noindent{\bf On Step 1}:

Let $r'\in [qr]$
and $\tilde r'\in[\tilde q\tilde r]$ with $|qr'|=|\tilde q\tilde r'|$ such that
$$|pr'|-|\tilde p\tilde r'|=\min\limits_{s\in[qr],
\tilde s\in[\tilde q\tilde r], |qs|=|\tilde q\tilde s|}\{|ps|-|\tilde p\tilde s|\}<0.$$
Then by Lemma 3.1 and 3.3 we can conclude that, similar to Lemma 1.2,
$$\angle qr'p=\angle \tilde q\tilde r'\tilde p\ \ and\ \ \angle rr'p=\angle \tilde r\tilde r'\tilde p.$$
On the other hand, since $|pr'|<|\tilde p\tilde r'|$, by Lemma 1.3 it is easy to see that
$$\angle \tilde q\tilde r'\tilde p<\tilde \angle qr'p\ \text{ or } \angle \tilde r\tilde r'\tilde p<\tilde \angle rr'p.$$
It therefore follows that
$$\angle qr'p<\tilde \angle qr'p \ \text{ or } \ \angle rr'p<\tilde \angle rr'p.$$

\vskip1mm

\noindent {\bf On Step 2}:

We can finish  Step 2 according to the following two cases.

\noindent Case 1: One of $[p_1q_1]$ and $[q_1r_1]$, say $[p_1q_1]$, is the minimal side of $\triangle p_1q_1r_1$.

We first notice that by (3.2) and Lemma 3.3 there is $t\in [q_1r_1]$ and
$\tilde t\in [\tilde q_1\tilde r_1]$ with $|\tilde t\tilde r_1|=|tr_1|$ such that $|p_1t|=|\tilde p_1\tilde t|$
and for any $t'\in [q_1t]^\circ$ and $\tilde t'\in [\tilde q_1\tilde t]$ with $|\tilde t'\tilde t|=|t't|$
$$|p_1t'|<|\tilde p_1\tilde t'|. \eqno{(3.5)}$$

If $t\neq r_1$, then $\text{peri}(p_1q_1t)<\text{peri}(p_1q_1r_1)$
and $\cos\angle p_1q_1t-\cos\tilde \angle p_1q_1t=\cos\angle p_1q_1r_1-\cos\tilde \angle p_1q_1r_1$.
I.e., it suffices to let $(p_2, q_2, r_2)=(p_1, q_1, t)$.

We now assume that $t=r_1$, and take the point $s\in [q_1r_1]$ with $|sr_1|=\frac12|p_1q_1|$.
It is clear that both $\text{peri}(p_1q_1s)$ and $\text{peri}(p_1r_1s)$ are less than $\text{peri}(p_1q_1r_1)$.
If $\angle p_1q_1s<\tilde \angle p_1q_1s$
or $\angle p_1sq_1<\tilde \angle p_1sq_1$,
then we can let $(p_2, q_2, r_2)\triangleq (p_1, q_1, s)$ or $(p_1, s, q_1)$ respectively
which satisfies (3.4) (due to (3.5), it is easy to see that
$\text{peri}(p_1q_1s)<\text{peri}(p_1q_1r_1)-\frac12|p_1q_1|$,
and if $(p_1q_1s)_{\min}<|p_1q_1|$ then $\text{peri}(p_1q_1s)<4|p_1q_1|$).
If $\angle p_1q_1s\geq\tilde \angle p_1q_1s$
and $\angle p_1sq_1\geq\tilde \angle p_1sq_1$, then we claim that $(p_2, q_2, r_2)=(p_1, s, r_1)$
which satisfies (3.3). First of all, note that $\angle p_1sr_1=\pi-\angle p_1sq_1$ and $\angle p_1q_1r_1=\angle p_1q_1s$.
Then by Lemma 3.4,
``$\angle p_1q_1s\geq\tilde \angle p_1q_1s$, $\angle p_1sq_1\geq\tilde \angle p_1sq_1$
and $\angle p_1q_1r_1<\tilde\angle p_1q_1r_1$'' implies that
$\angle p_1sr_1<\tilde \angle p_1sr_1$, and thus one can check that
$$\frac{\cos\angle p_1sr_1-\cos\tilde\angle p_1sr_1}{\cos\angle p_1q_1r_1-\cos\tilde \angle p_1q_1r_1}\ge
\frac{\cos(\pi-\tilde \angle p_1sq_1)-\cos\tilde\angle p_1sr_1}{\cos\tilde\angle p_1q_1s-\cos\tilde \angle p_1q_1r_1}
=\frac{|p_1q_1|\cdot|q_1r_1|}{|p_1s|\cdot|sr_1|}\stackrel{\text{by (3.5)}}{>}1.$$

\vskip1mm

\noindent  Case 2: $[p_1r_1]$ is the strictly minimal side of $\triangle p_1q_1r_1$.

Without loss of generality, we assume that $|p_1r_1|<|p_1q_1|\leq |q_1r_1|$.
By the same reason as in Case 1, we can assume that (3.5) holds for
all $t'\in [q_1r_1]^\circ$ and $\tilde t'\in [\tilde q_1\tilde r_1]$ with $|\tilde t'\tilde r_1|=|t'r_1|$.
Then we can conclude that there is $s\in [q_1r_1]$ such that either $|sr_1|=|p_1r_1|$ and $|p_1r_1|\leq |p_1s|<\sqrt2|p_1r_1|$,
or $|sr_1|>|p_1r_1|$ and $|p_1s|=|p_1r_1|$.
It follows that $\text{peri}(p_1q_1s)<\text{peri}(p_1q_1r_1)-\frac12|p_1r_1|$,
and $\text{peri}(p_1sq_1)<4.5|p_1r_1|$ if $|q_1s|<|p_1r_1|$.
Moreover, similar to the proof for Case 1, we can let $(p_2, q_2, r_2)=(p_1, q_1, s)$ or $(p_1, s, q_1)$
which satisfies (3.4), or let $(p_2, q_2, r_2)=(p_1, s, r_1)$ which satisfies
(3.3).
\hfill$\Box$

\vskip2mm

\noindent {\bf Remark 3.5} Let $\{\triangle p_iq_ir_i\}_{i=1}^\infty$ be
triangles in the right above proof, and
define $d_{ij}\triangleq\max\{|p_j\{p_i,q_i,r_i\}|,|q_j\{p_i,q_i,r_i\}|, |r_j\{p_i,q_i,r_i\}|\}$ for $j>i\geq 1$.
We can show that $\{d_{1j}\}_{j=2}^\infty$ have an upper bound, so
$\{\triangle p_iq_ir_i\}_{i=1}^\infty$ lie in a compact domain of $X$ even when $X$ is not compact.
Note that by the choice of $\triangle p_iq_ir_i$,
it is possible that
$$d_{i(i+1)}\leq\text{peri}(p_iq_ir_i)-\text{peri}(p_{i+1}q_{i+1}r_{i+1});\eqno{(3.6)}$$
and thus it is clear that if (3.6) holds for all $i$ with $1\leq i\leq j-1$, then
$$d_{1j}\leq \sum\limits_{i=1}^{j-1}d_{i(i+1)}
\leq\sum\limits_{i=1}^{j-1}\left(\text{peri}(p_iq_ir_i)-\text{peri}(p_{i+1}q_{i+1}r_{i+1})\right)
<\text{peri}(p_1q_1r_1).$$
Unfortunately, (3.6) might not be true,
and even there is no constant $c\in(0,1)$ such that
$c\cdot d_{i(i+1)}\leq\text{peri}(p_iq_ir_i)-\text{peri}(p_{i+1}q_{i+1}r_{i+1})$ for all $i$.
However, by analyzing the process of locating $\triangle p_iq_ir_i$ in Step 2 and 3 case by case,
we can conclude that, for all $j$, there is a constant $c\in(0,1)$ and $1=i_0<i_1<\cdots<i_k=j$ such that
$$d_{1j}\leq\sum\limits_{l=0}^{k-1}d_{i_{l}i_{l+1}}
\leq\frac{1}{c}\cdot\sum\limits_{l=0}^{k-1}\left(\text{peri}(p_{i_l}q_{i_l}r_{i_l})-\text{peri}(p_{i_{l+1}}q_{i_{l+1}}r_{i_{l+1}})\right)
<\frac{1}{c}\cdot\text{peri}(p_1q_1r_1)
$$
(where it might hold that $c\cdot d_{i_{l}i_{l+1}}>\text{peri}(p_{i_l}q_{i_l}r_{i_l})-\text{peri}(p_{i_{l+1}}q_{i_{l+1}}r_{i_{l+1}})$ for some $l$).


\noindent School of Mathematical Sciences (and Lab. math. Com.
Sys.), Beijing Normal University, Beijing, 100875
P.R.C.\\
e-mail: wyusheng$@$bnu.edu.cn

\end{document}